\documentclass{article}%
\usepackage{amssymb}
\usepackage{amsfonts}
\usepackage{amsmath}%
\usepackage{graphicx}
\providecommand{\U}[1]{\protect\rule{.1in}{.1in}}
\newtheorem{theorem}{Theorem}
\newtheorem{acknowledgement}[theorem]{Acknowledgement}

\begin{document}

\title{The High Precision Numerical Calculation of Stieltjes Constants. Simple and
Fast Algorithm}
\author{Krzysztof Ma\'{s}lanka\\Polish Academy of Sciences\\Institute for the History of Science\\Nowy \'{S}wiat 72, 00-330 Warsaw, Poland\\e-mail krzysiek2357@gmail.com
\and Andrzej Kole\.{z}y\'{n}ski\\University of Science and Technology\\Faculty of Materials Science and Ceramics\\Mickiewicza 30, 30-059 Cracow, Poland\\e-mail kolezyn@agh.edu.pl}
\maketitle

\begin{abstract}
We present a simple but efficient method of calculating Stieltjes constants at
a very high level of precision, up to about 80000 significant digits. This
method is based on the hypergeometric-like expansion for the Riemann zeta
function presented by one of the authors in 1997 \cite{Maslanka 1}. The
crucial ingredient in this method is a sequence of high-precision numerical
values of the Riemann zeta function computed in equally spaced real arguments,
i.e. $\zeta(1+\varepsilon),\zeta(1+2\varepsilon),\zeta(1+3\varepsilon),...$
where $\varepsilon$ is some real parameter. (Practical choice of $\varepsilon
$\ is described in the main text.) Such values of zeta may be readily obtained
using the PARI/GP program, which is especially suitable for this.


\end{abstract}

\section{Introduction: the Riemann $\zeta$ function}

Fundamental formulas in number theory are seldom numerically efficient.
Although deep and absolutely precise, they may even hide the most important
features of involved quantities. As a prominent example we consider the
celebrated zeta function $\zeta(s)$ discovered by Euler in 1737 and published
in 1744 \cite{Euler} as a function of real variable and meticulously
investigated by Riemann in the complex domain in his famous memoir submitted
in 1859 to the Prussian Academy \cite{Riemann}:%
\begin{equation}
\zeta(s)=%
{\displaystyle\sum\limits_{n=1}^{\infty}}
\frac{1}{n^{s}}\qquad\operatorname{Re}(s)>1\label{Zeta}%
\end{equation}
This is a special case of more general class of functions called Dirichlet
series. It is divergent in the most interesting area of the complex plane,
i.e., in the so called critical strip $0\leq\operatorname{Re}s\leq1$ where all
complex zeros of zeta lie. However, as was shown by Riemann, the definition
(\ref{Zeta}) does contain information about the zeta function on the entire
complex plane but the process of analytic continuation must be used in order
to reveal global behavior of this function. There is no universal procedure
how to achieve this in practice and usually various ingenious tricks are
required. For example, considering simply alternating version of (\ref{Zeta})
leads to another Dirichlet series which is convergent for $\operatorname{Re}%
s>0$ (except $s=1$), i.e.\ also inside the critical strip:%
\[
\zeta(s)=\frac{1}{1-2^{1-s}}%
{\displaystyle\sum\limits_{n=1}^{\infty}}
\frac{(-1)^{n}}{n^{s}}\qquad\operatorname{Re}(s)>0,\quad s\neq1
\]
However, in order to obtain globally convergent representation for $\zeta
$\ one has to use more sophisticated techniques. We shall describe such an
approach below.%
\begin{figure}[ptb]%
\centering
\includegraphics[height=0.5\textheight]{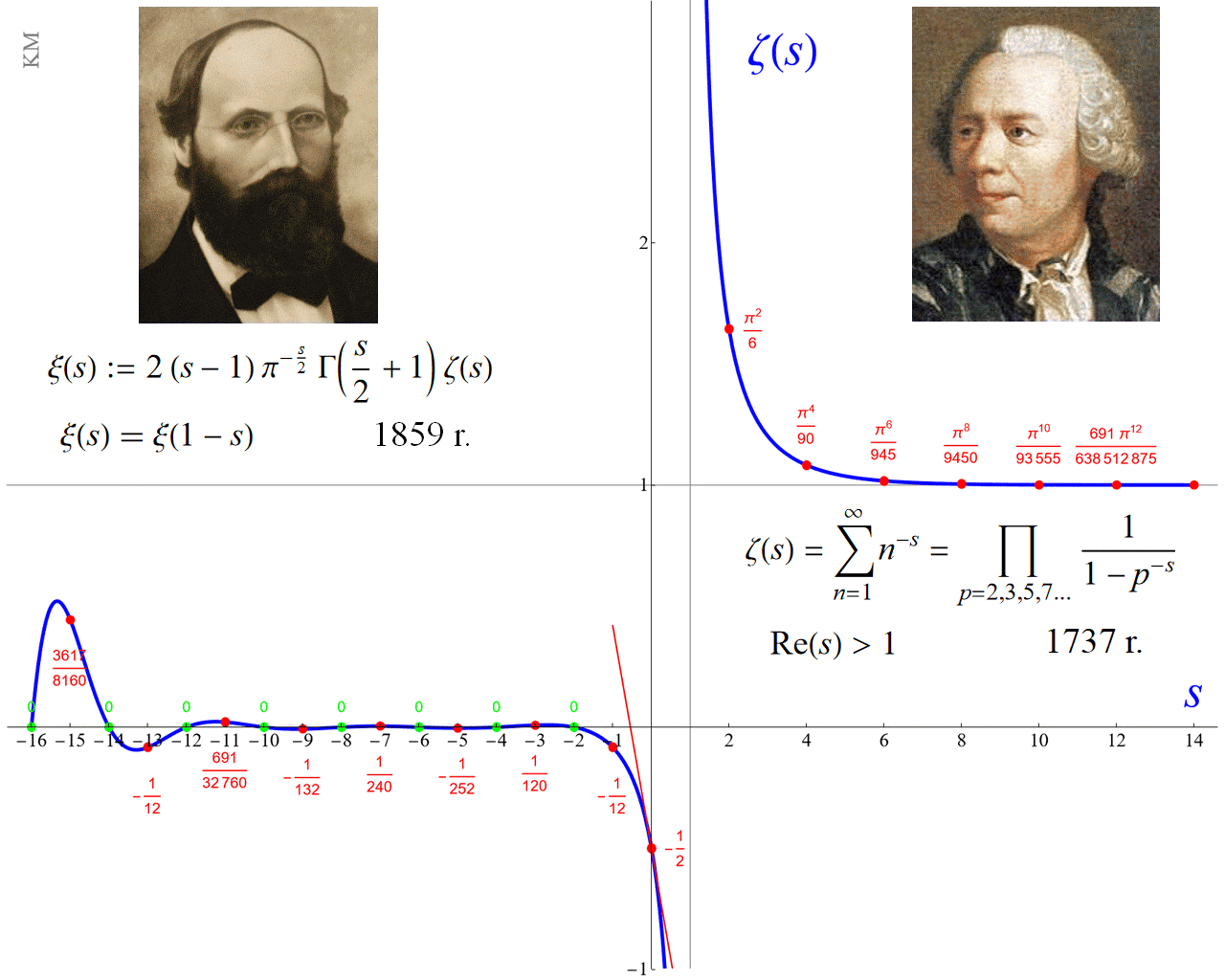}%
\caption{Plot of the zeta function for real variable (blue curve). Euler
discovered the zeta function in 1737 and found its deep connection with prime
numbers \cite{Euler} (see box on the right). But it was Riemann who in 1859
rigorously proved certain fundamental equation for it and made its analytical
continuation to the entire complex plane, except for a single pole for $s=1$
\cite{Riemann} (box on the left). Values of zeta for $s=2n$, $n=1,2,...$ were
found by Euler in closed form (red dots). $\zeta(-2n)=0$ are so called trivial
zeros (green dots).}%
\end{figure}

The Riemann zeta function contains the (heavily encoded) puzzle of the
distribution of prime numbers. According to the famous saying of Paul
Erd\"{o}s (1913-1996) -- the solution to this puzzle may appear only "in
millions of years, but even then it will not be complete, because in this case
we are facing Infinity". We know, however, that this secret lies in the
distribution of the zeros of the zeta function, i.e. the roots of the "simple"
equation $\zeta(s)=0$, on the complex plane. In 1859 Riemann hypothesized that
\textit{all} these roots (except for the so-called trivial ones) lie precisely
on the line $\operatorname{Re}s=\frac{1}{2}$.%

\begin{figure}[ptb]%
\centering
\includegraphics[height=0.6\textheight]{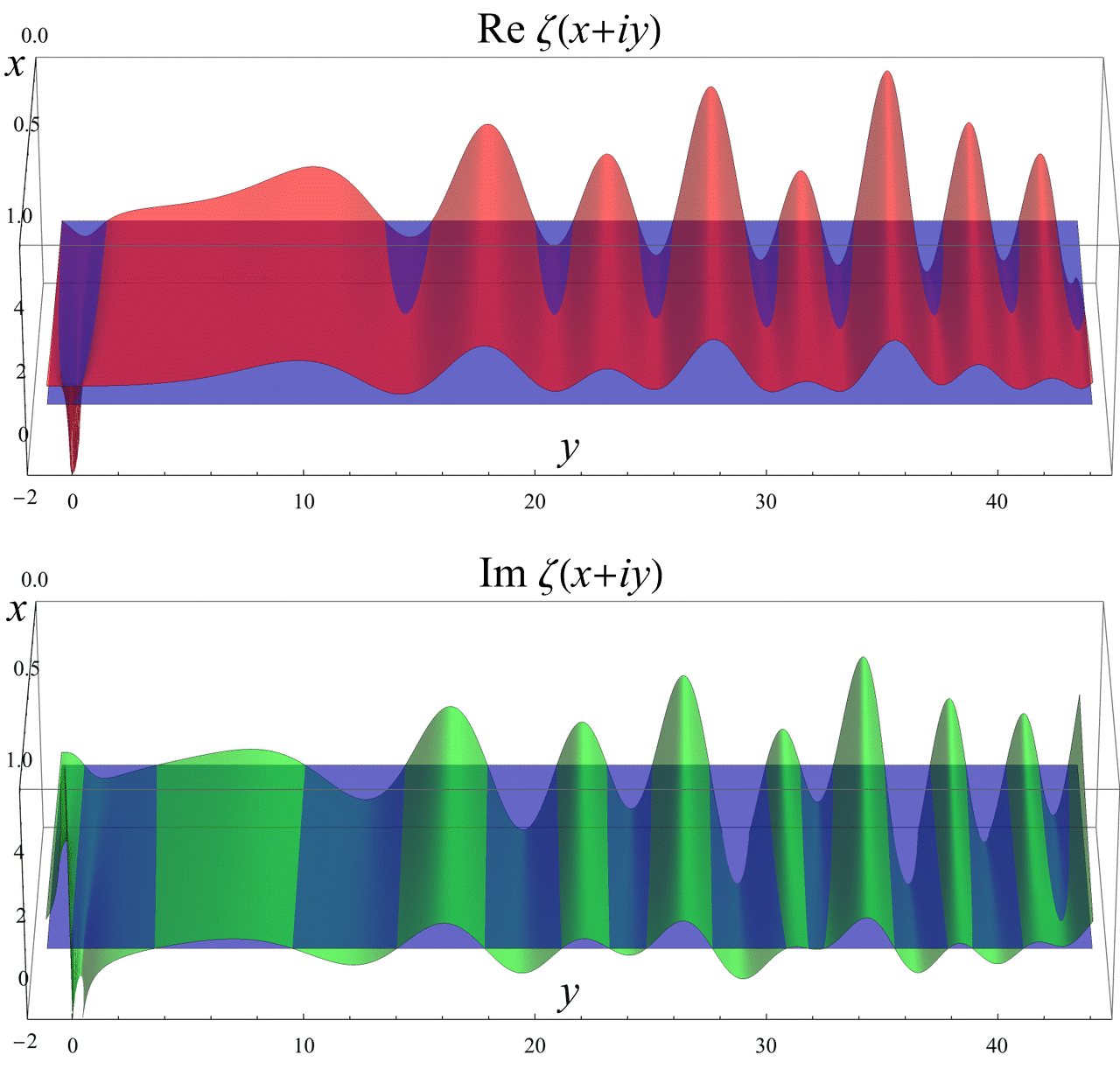}%
\caption{The zeta function shows its essence and its true meaning only in the
complex domain, and we owe knowledge about it to Riemann. The upper graph is
the real part of the zeta function $\zeta(s)$, the lower graph is its
imaginary part in the complex domain. The blue plane is the plane of the
complex variable $s$.}%
\end{figure}
\begin{figure}[ptb]%
\centering
\includegraphics[height=0.3\textheight]{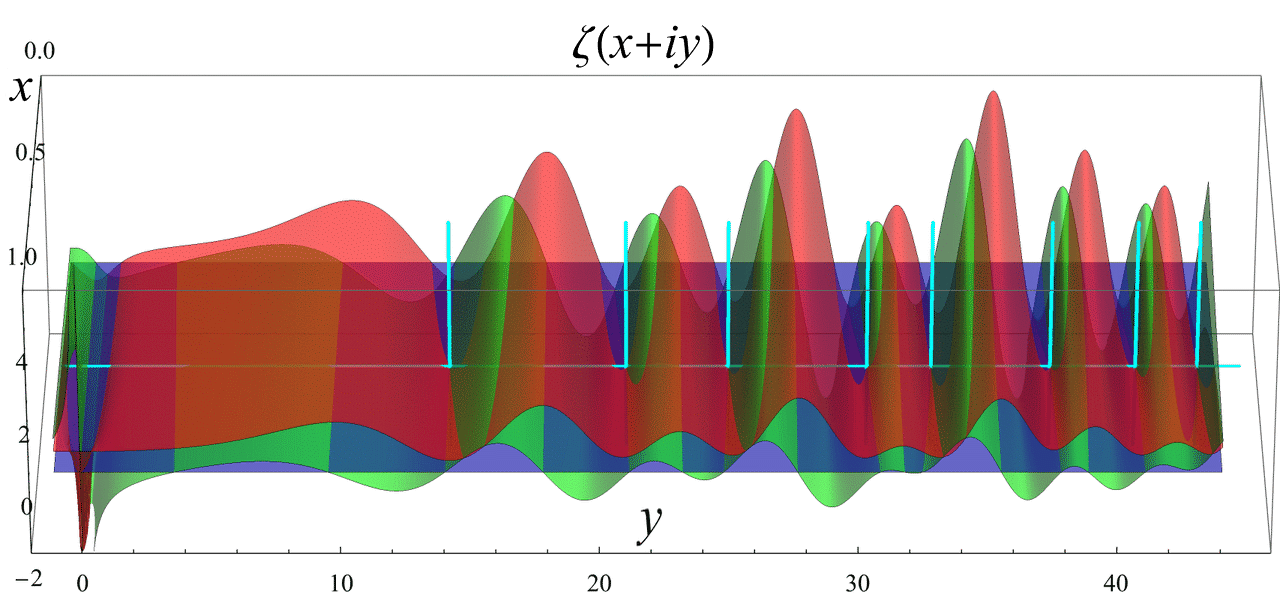}%
\caption{Both surfaces shown in Fig. 2 intersects the plane of the complex
variable $s$ along certain irregular curves.After overlapping these surfaces,
it turns out that these curves intersect themselves at certain points -- these
are the complex zeros of the zeta function (indicated by vertical blue lines).
The Riemann hypothesis says that \textit{all} these zeros are placed exactly
on the line $\operatorname{Re}s=1/2$.}%
\end{figure}
Despite the passage of more than a century and a half and the persistent
efforts of many top-class mathematical talents, the Riemann hypothesis remains
unsettled. We simply do not know whether it is true or false. (Some think that
it is undecidable.) Computer experiments based on billions of numerically
calculated complex roots seem to confirm it. However, exact proof remains, so
far, beyond the reach of mathematicians. It seems no one has even had a good
idea of how to attack this problem so far. Some have suggested that some "new
math" is needed for this, but this view is too vague to be of any practical help.

\section{Stieltjes constants}

The Stieltjes constants are closely related to the Riemann zeta function, and
since this function is extremely important in analytical number theory, these
constants are equally important.

Formulas for the Stieltjes constants may serve as another example of strict
and deep but numerically inefficient formulas. These constants are essentially
coefficients of the Laurent series expansion of the zeta function around its
only simple pole at $s=1$:%
\begin{equation}
\zeta(s)=\frac{1}{s-1}+%
{\displaystyle\sum\limits_{n=0}^{\infty}}
\frac{\left(  -1\right)  ^{n}}{n!}\gamma_{n}\left(  s-1\right)  ^{n}%
\label{ZetaExpansion}%
\end{equation}
Primary definition of these fundamental constants was found by Thomas Jan
Stieltjes and presented in a letter to his close friend and collaborator
Charles Hermite dated June 23, 1885 \cite{Hermite}:%
\begin{equation}
\gamma_{n}=\underset{m\rightarrow\infty}{\lim}\left(
{\displaystyle\sum\limits_{k=1}^{m}}
\frac{\left(  \ln k\right)  ^{n}}{k}-\frac{\left(  \ln m\right)  ^{n+1}}%
{n+1}\right) \label{gamma}%
\end{equation}
When $n=0$ the numerator in the first summand in (\ref{gamma}) is formally
$0^{0}$ which is taken to be $1$. In this case, (\ref{gamma}) reduces simply
to the well-known Euler-Mascheroni constant%
\[
\gamma_{0}=\underset{m\rightarrow\infty}{\lim}\left(
{\displaystyle\sum\limits_{k=1}^{m}}
\frac{1}{k}-\ln m\right)
\]
which, roughly speaking, measures the rate of divergence of the harmonic series.

Effective numerical computing of the constants $\gamma_{n}$ is quite a
challenge because the formulas (\ref{gamma}) are extremely slowly convergent.
Even for $n=0$, in order to obtain just $10$ accurate digits one has to sum up
exactly $12366$ terms whereas in order to obtain $10000$ digits (which is
indeed required in some applications) one would have to sum up unrealistically
large number of terms: nearly $5\cdot10^{4342}$ which is of course far beyond
capabilities of the present day computers. For $n>0$ the situation is still
worse. Therefore we have to seek for other faster algorithms.

Due to the terribly slow convergence mentioned above, the progress in
calculating the numerical values of Stieltjes constants values was very slow.
In his letter to Hermite Stieltjes himself gave just two very inaccurate
values for these constants$\gamma_{n}$ (except for the then well-known
Euler-Mascheroni constant $\gamma_{0}$):%
\begin{align*}
\gamma_{1}  & =-0.072815(520)...\\
\gamma_{2}  & =-0.004(7)...
\end{align*}
(Here and below digits in brackets are incorrect.) Two years later, in 1887,
Jensen \cite{Jensen} gave eight values with nine significant digits:%
\begin{align*}
\gamma_{1}  & =+0.072815845...\\
\gamma_{2}  & =-0.004845182...\\
\gamma_{3}  & =-0.000342306...\\
\gamma_{4}  & =+0.0000968(89)...\\
\gamma_{5}  & =-0.000006611...\\
\gamma_{6}  & =-0.000000332...\\
\gamma_{7}  & =+0.000000105...\\
\gamma_{8}  & =-0.000000009...
\end{align*}

Certain hope is in using integral representations of the Stieltjes constants.
There are at least three such integrals:

-- by directly applying Cauchy integral formula for derivatives to the Riemann
zeta function we get:%

\begin{equation}
\gamma_{n}=\frac{(-1)^{n}n!}{2\pi}\int_{0}^{2\pi}e^{-nit}\zeta\left(
e^{it}+1\right)  dt\label{Cauchy}%
\end{equation}

-- by Franel, 1895, \cite{Franel}%

\begin{equation}
\gamma_{n}=\frac{1}{2}\delta_{n,0}+\frac{1}{i}\int_{0}^{\infty}\frac
{dt}{e^{2\pi t}-1}\left[  \frac{\left(  \ln\left(  1-it\right)  \right)  ^{n}%
}{1-it}-\frac{\left(  \ln\left(  1+it\right)  \right)  ^{n}}{1+it}\right]
\label{Franel}%
\end{equation}

-- by Blagouchine \cite{Blagouchine}%

\begin{equation}
\gamma_{n}=\frac{\pi}{2(n+1)}\int_{-\infty}^{\infty}\frac{\left(  \ln\left(
\frac{1}{2}\pm it\right)  \right)  ^{n+1}}{\left(  \cosh\pi t\right)  ^{2}%
}dt\label{Blagouchine}%
\end{equation}

Using these integral representations one can, for exaple with the help of the
procedure \textbf{NIntegrate} which is built in Wolfram \textit{Mathematica},
calculate $\gamma_{n}$ up to $n=1000$ with precision of several hundred
significant digits in a reasonable computer time. However, increasing $n$
and/or the working precision parameter in \textbf{NIntegrate} produces error
message in the \textit{Mathematica} output. This can easily be understood when
looking at the behavior of the integrand of (\ref{Cauchy}) which for growing
$n$ contains more large oscillations.%

\begin{figure}[ptb]%
\centering
\includegraphics[height=0.5\textheight]{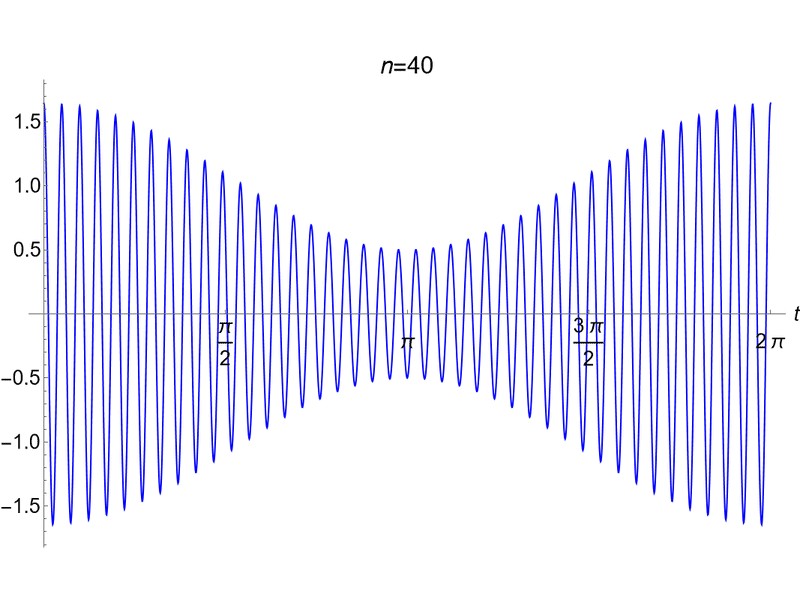}%
\caption{Real part of the integrand (\ref{Cauchy}) which contributes to the
value of $\gamma_{40} $. (Integrating the imaginary part which is
antisymmetric with respect to $t=\pi$ gives zero.) Number of oscillations
grows as $n$. Therefore for large $n$ the numerical integration procedure
cannot properly estimate this integral.}%
\end{figure}
%

\begin{figure}[ptb]%
\centering
\includegraphics[height=0.5\textheight]{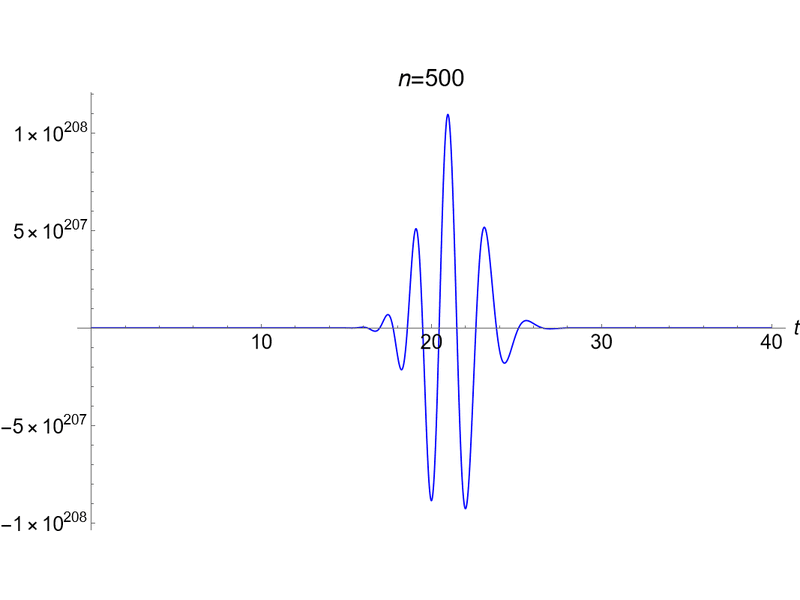}%
\caption{Real part of the integrand in (\ref{Blagouchine}). The number of
oscillations grows with $n$ but, contrary to the case (\ref{Franel}), their
amplitude also increases very quickly with $n$.}%
\end{figure}

There are also several series representations of Stieltjes constants, e.g.
such as this given by I. Blagouchine \cite{Blagouchine}%

\begin{equation}
\gamma_{n}=-\frac{1}{n+1}\sum_{k=0}^{\infty}\frac{1}{k+1}\sum_{j=0}^{\infty
}(-1)^{j}\binom{k}{j}\log^{n+1}\left(  j+1\right) \label{Blagouchine 1}%
\end{equation}
and another one found by M. Coffey \cite{Coffey} (Corollary 13, with misprint)%
\begin{equation}
\gamma_{n}=-\frac{2}{3}n!\sum_{i=1}^{n}\frac{B_{n-i+1}c^{n-i}}{i!(n-i+1)!}%
h_{i}-\frac{2}{3c(n+1)}h_{n+1}-\frac{B_{n+1}}{n+1}c^{n+1}\label{Coffey}%
\end{equation}
where%
\begin{align*}
c  & \equiv\ln2\\
h_{i}  & \equiv\sum_{k=1}^{\infty}3^{-k}\sum_{j=1}^{k}\frac{(-1)^{j}2^{j}%
}{j+1}\binom{k}{j}\ln^{i}(j+1)
\end{align*}
Unfortunately, both (\ref{Blagouchine 1}) and (\ref{Coffey}) are also very
slowly convergent and pretty useless in numerical investigations -- contrary
to what Coffey claims:\ "The expression may be attractive for some
computational applications because it exhibits even faster convergence" (see
\cite{Coffey}, p. 23).

Significant progress took place in 1984-1985 with the work of Ainsworth and
Howell \cite{Ainsworth Howell}\ who got a grant from NASA and probably used a
computer. (It could be an analog machine, but they did not disclose the
technical details of their calculations.) They used another integral
representation of the Stieltjes constants and with the help of the Gauss
numerical integration formula tabulated $200$ initial $\gamma_{n}$ with just
$10$ significant digits each. They also calculated a few selected values of
$\gamma_{n}$ for larger $n=500,1000,1500,2000$. In the latter cases, some of
their digits are incorrect.

In 1992 Keiper\footnote{Jerry B. Keiper (1953-1995) worked for Wolfram
Research and was an active contributor to \textit{Mathematica}. He developed,
among others many effective algorithms for numerical computation of special
functions. He died tragically returning from work on his bike, hit by a car.}
published an effective algorithm for calculating Stieltjes constants. Keiper's
algorithm was later implemented in \textit{Mathematica \cite{Keiper}}.
(However, no technical details about this algorithm can be found in
\textit{Mathematica} documentation except for a concise statement that it
"uses Keiper's algorithm based on numerical quadrature of an integral
representation of the zeta function and alternating series summation using
Bernoulli numbers".)

An efficient but rather complicated method based on Newton-Cotes quadrature
has been proposed by Kreminski in 2003 \cite{Kreminski}. This was a real
achievement since\ Kreminski computed $\gamma_{n}$ up to $n=3000$ with several
thousand digits and was able to observe certain interesting structures in the
distribution of $\gamma_{n}$.

Quite recently (2013) Johansson presented particularly efficient method
\cite{Johansson}. He calculated a new, impressive, record-breaking value of
$\gamma_{n}$ for $n=100000$. Later (2018), in collaboration with Blagouchine,
Johansson reached the next record values: $n=10^{10},10^{15}$ and $10^{100}$.

In the present paper yet another method of computing Stieltjes constants will
be described which, we believe, is perhaps not as efficient as Johansson's
approach, yet it is by far more simple and it may be easily and quickly used
in practical calculations for obtaining $\gamma_{n}$ up to $n\sim30000$ with
precision $\sim80000$ significant digits.

\section{Riemann zeta representation}

In 1997 it was shown by one of the authors of the present paper \cite{Maslanka
1}, \cite{Maslanka 2} that the Riemann zeta function may be expressed as%
\begin{align}
\zeta(s)  & =\frac{1}{s-1}\left[  A_{0}+\left(  1-\frac{s}{2}\right)
A_{1}+\left(  1-\frac{s}{2}\right)  \left(  2-\frac{s}{2}\right)  \frac{A_{2}%
}{2!}+...\right]  =\label{Hyper}\\
& =\frac{1}{s-1}%
{\displaystyle\sum\limits_{k=0}^{\infty}}
\frac{A_{k}}{k!}%
{\displaystyle\prod\limits_{i=1}^{k}}
\left(  i-\frac{s}{2}\right)  =\label{Hyper2}\\
& =\frac{1}{s-1}%
{\displaystyle\sum\limits_{k=0}^{\infty}}
\frac{\Gamma\left(  k+1-\frac{s}{2}\right)  }{\Gamma\left(  1-\frac{s}%
{2}\right)  }\frac{A_{k}}{k!}\label{Hyper3}\\
& =\frac{1}{s-1}%
{\displaystyle\sum\limits_{k=0}^{\infty}}
\left(  1-\frac{s}{2}\right)  _{k}\frac{A_{k}}{k!}\qquad s\in%
\mathbb{C}
\backslash\{1\}\label{Hyper4}%
\end{align}
where%
\begin{align}
A_{k}  & =%
{\displaystyle\sum\limits_{j=0}^{k}}
\left(  -1\right)  ^{j}\binom{k}{j}(2j+1)\zeta(2j+2)=\label{Ak}\\
& =\frac{1}{2}%
{\displaystyle\sum\limits_{j=0}^{k}}
\binom{k}{j}(2j+1)\frac{\left(  2\pi\right)  ^{2j+2}B_{2j+2}}{\left(
2j+2\right)  !}%
\end{align}
$\left(  x\right)  _{k}$ is the Pochhammer symbol and $B_{n}$ is the
$n^{\text{th}}$ Bernoulli number \cite{Abramowitz}.

The main idea behind this approach is to remove the single pole of the zeta
function multiplying it by $s-1$ and then to fix values of this entire
function in an infinite number of \textit{equally spaced} real points which
corresponds simply to interpolation with nodes. Note that in (\ref{Hyper}%
)--(\ref{Hyper4}) these points are precisely the points in which, as shown by
Euler, zeta values are known exactly, i.e. $s=2,4,6,...$ Indeed, series
(\ref{Hyper}) truncates in these points and gives appropriate exact values.

It may be shown that this representation is globally convergent. Real
coefficients $A_{k}$ expressed as an alternating binomial sum are in fact
combinations of Bernoulli numbers and even powers of $\pi$. On the other hand,
$\left(  -1\right)  ^{k}A_{k}$ are simply consecutive finite step derivatives
of some entire function, namely $(2s+1)\zeta(2s+2)$ that involves these points
in which, as Euler had shown, zeta is explicitly known.

Considered as a sort of polynomial interpolation with fixed nodes the
expansion (\ref{Hyper}) might appear trivial. However, this is not the case
since many "simple" functions, e.g. Lorentz function $1/(1+x^{2})$, exhibit
nasty phenomenon known as the Runge effect: oscillations between fixed nodes
growing when number of terms in the series increases. This behavior may be
cured using \textit{unequally} spaced nodes, so called Chebyshev nodes, but
this in turn spoils the very idea of (\ref{Hyper}) which leads in a natural
way to Pochhammer symbols. From this point of view the global validity of
(\ref{Hyper}) is equivalent to the following simple statement: the regularized
Riemann zeta function $(s-1)\zeta(s)$\ does not exhibit Runge phenomenon.

The original proof of (\ref{Hyper}) contained a gap \cite{Maslanka 1}.
Rigorous proof has been given by B\'{a}ez-Duarte\footnote{The prominent
Venezuelan mathematician Luis B\'{a}ez-Duarte (1938-2018), educated in the
USA, Massachusetts Institute of Technology, and working at the Instituto
Venezolano de Investigaciones Cient\'{\i}ficas (IVIC) in Caracas, was a close
friend and collaborator of one of the authors (K.M.). Although they had never
met in person, from 2003 until Luis' death, they conducted lively
correspondence, mainly on mathematical topics, but also on general topics
related to literature, history, politics, etc.} in 2003 \cite{Baez 1} who also
presented certain simple and esthetic\ criterion for the Riemann Hypothesis
based on expansion (\ref{Hyper}) \cite{Baez 2}. Another very short and
particularly elegant proof of (\ref{Hyper}) using Carlson theorem has been
given by Flajolet and Vepstas in 2007 \cite{Flajolet}. Later a whole class of
similar zeta representations has been published \cite{Maslanka 2}.

Coefficients $A_{k}$ tend to zero sufficiently fast which is crucial to assure
the global convergence of the series (\ref{Hyper}). However, their detailed
behavior with growing $k$ is quite striking as can be seen on a logarithmic
plot with the $k$-axis rescaled as $\sqrt[3]{k}$. More precisely: they exhibit
curious and unexpected oscillatory behavior with both amplitude and frequency
decreasing when $k$ tends to infinity (see Fig. 6).%
\begin{figure}[ptb]%
\centering
\includegraphics[height=0.4\textheight]{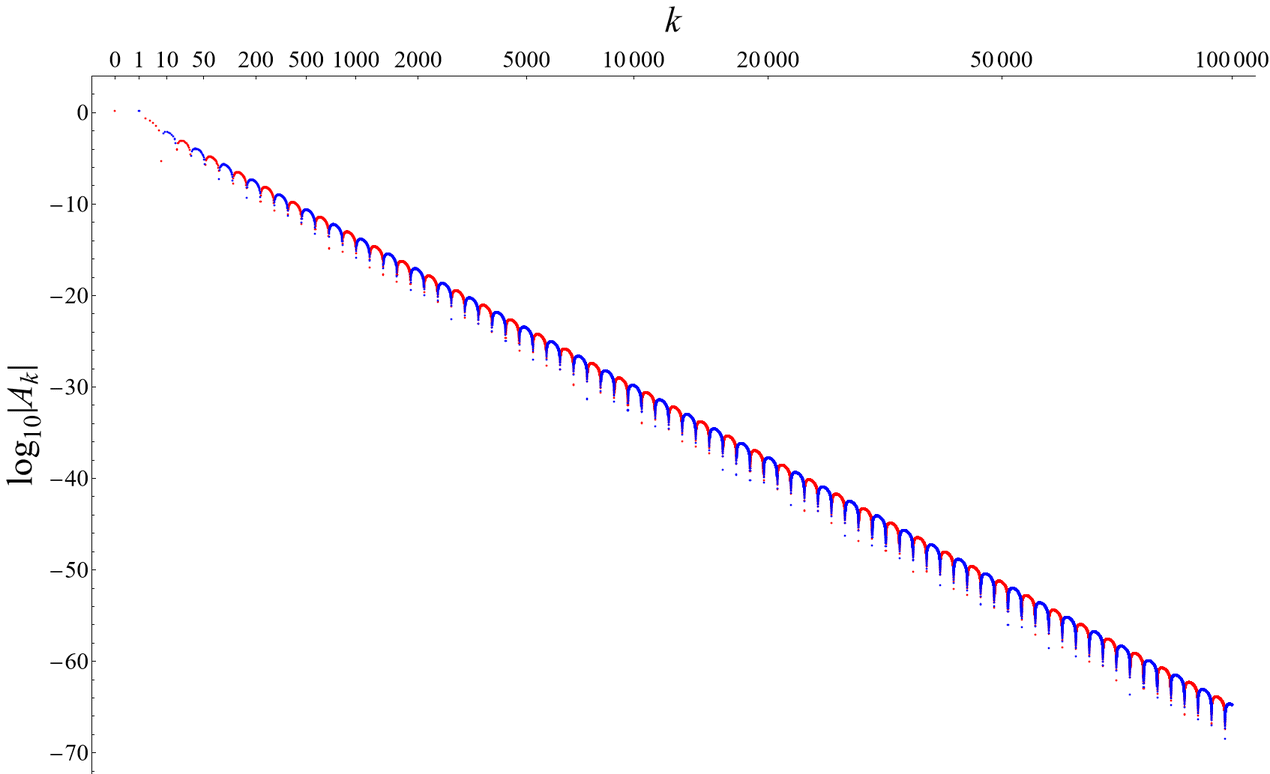}%
\caption{Curious behavior of the coefficients $A_{k}$. given by (\ref{Ak}).
There are unexpected oscillations with slowly diminishing frequency (roughly
as $k^{-2/3}$) and nearly exponentially diminishing amplitude. Note that the
$k$-axis is scaled as $\sqrt[3]{k}$. Red points correspond to positive values
of $A_{k}$ and blue to negative ones.}%
\end{figure}
This peculiar behavior "cries for explanation" as stated in \cite{Flajolet}
(p. 2). Using the saddle point method one can show \cite{Maslanka 3} that for
$k$ tending to infinity the following asymptotics holds:%
\begin{equation}
A_{k}\sim\frac{4\pi^{3/2}}{\sqrt{3\kappa}}\exp\left(  -\frac{3}{2}\kappa
+\frac{\pi^{2}}{4\kappa}\right)  \cos\left(  \frac{4\pi}{3}-\frac{3\sqrt{3}%
}{2}\kappa-\frac{\sqrt{3}\pi^{2}}{4\kappa}\right) \label{Asymptotic 1}%
\end{equation}

where%
\[
\kappa\equiv\pi^{3/2}\sqrt[3]{k}%
\]%
\begin{figure}[ptb]%
\centering
\includegraphics[height=0.4\textheight]{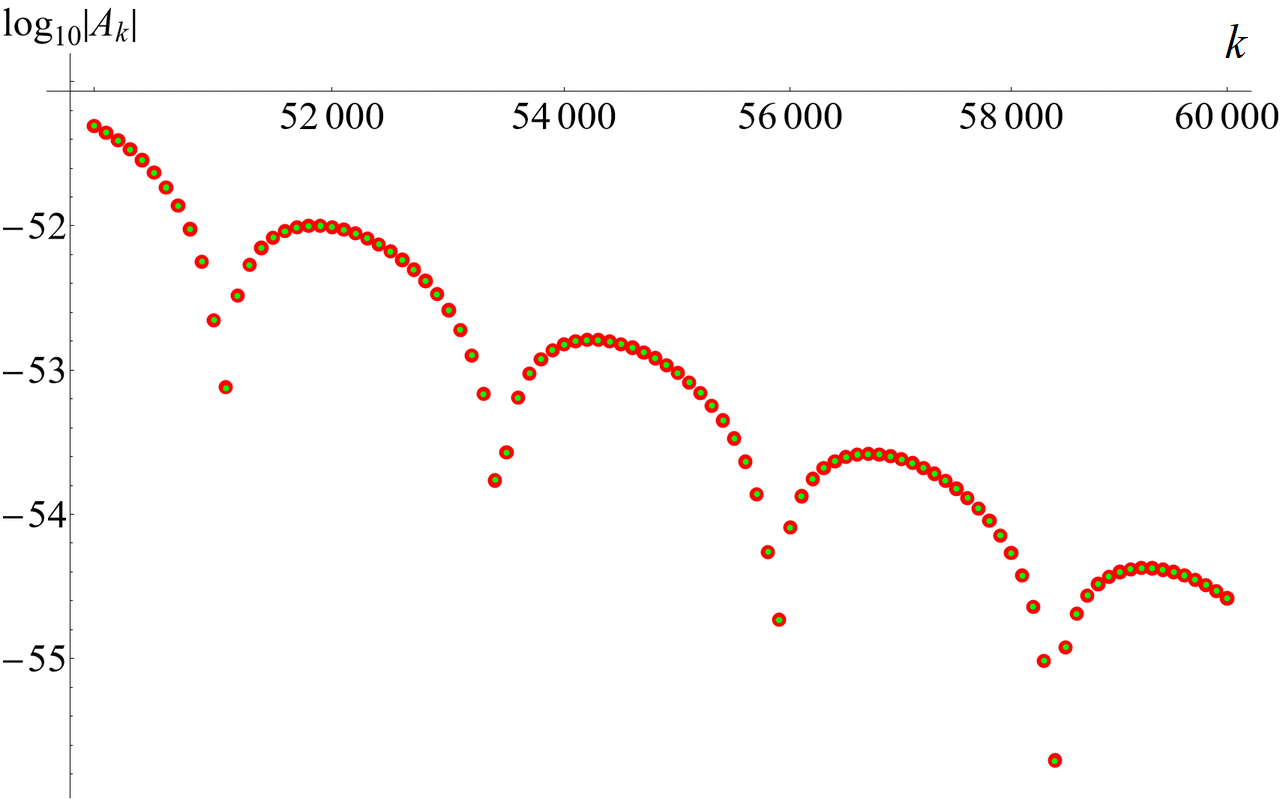}%
\caption{Asymptotic formula (\ref{Asymptotic 1}) works pretty well: red dots
represent exact values of $A_{k}$ as given by (\ref{Ak}) whereas smaller green
dots are calculated from (\ref{Asymptotic 1}).}%
\end{figure}

Coefficients $A_{k}$ obey also certain simple algebraic identities which stem
directly from trivial zeros of zeta and from the fact that $\zeta(0)=-\frac
{1}{2}$. Indeed, substituting in (\ref{Hyper}) $s=0,-2,-4,-6,...$ and making
use of elementary properties of the Euler gamma function we successively get:%
\begin{align}%
{\displaystyle\sum\limits_{k=0}^{\infty}}
A_{k}  & =\frac{1}{2}\label{relations}\\%
{\displaystyle\sum\limits_{k=0}^{\infty}}
(k+1)A_{k}  & =0\nonumber\\%
{\displaystyle\sum\limits_{k=0}^{\infty}}
(k+1)(k+2)A_{k}  & =0\nonumber\\%
{\displaystyle\sum\limits_{k=0}^{\infty}}
(k+1)(k+2)(k+3)A_{k}  & =0\nonumber\\
& ...\nonumber
\end{align}
After some simple manipulations we finally get:%
\begin{equation}%
{\displaystyle\sum\limits_{k=0}^{\infty}}
k^{n}A_{k}=\frac{(-1)^{n}}{2}\qquad n=0,1,2,...\label{Identities}%
\end{equation}
with the convention $k^{n}=1$ when $k=n=0$. Unfortunately, due to slow
convergence of (\ref{Identities}) when $n$ is large, these identities cannot
be effectively used to calculate $A_{k}$. Another interesting identity follows
from $\zeta^{\prime}(0)=-\frac{1}{2}\log(2\pi)$:%
\begin{equation}%
{\displaystyle\sum\limits_{k=0}^{\infty}}
A_{k}H_{k}=1-\log(2\pi)\label{identity}%
\end{equation}
where $H_{k}\equiv%
{\textstyle\sum\limits_{i=1}^{k}}
\frac{1}{i}$ is the $k^{\text{th}}$ harmonic number.

\section{Algorithm for calculating Stieltjes constants}

The particular choice of nodes in $s=2,4,6,...$ in the expansion
(\ref{Hyper}), albeit the most natural, is by no means the only one. One only
requires that the prescribed points be strictly equally spaced. For the
purpose of present calculations we choose the following sequence of points:%
\[
1,1+\varepsilon,1+2\varepsilon,1+3\varepsilon,...
\]
where $\varepsilon$ is certain real, not necessarily small number.

More precisely, define certain entire function $f$ as:%
\begin{equation}
f(s):=\left\{
\genfrac{}{}{0pt}{}{\zeta(s)-\frac{1}{s-1}\qquad s\neq1}{\gamma\qquad
\qquad\qquad s=1}%
\right. \label{regularized}%
\end{equation}
where $\gamma$ is the Euler constants which stems from the appropriate limit.
Then, instead of (\ref{Hyper3}), we have%
\[
f(s)=%
{\displaystyle\sum\limits_{k=0}^{\infty}}
\frac{\Gamma\left(  k-\frac{s-1}{\varepsilon}\right)  }{\Gamma\left(
-\frac{s-1}{\varepsilon}\right)  }\frac{\alpha_{k}}{k!}%
\]
with%
\begin{equation}
\alpha_{k}=%
{\displaystyle\sum\limits_{j=0}^{k}}
\left(  -1\right)  ^{j}\binom{k}{j}f(1+j\varepsilon)\label{Ak general}%
\end{equation}
Note that coefficients $\alpha_{k}$\ depend on $\varepsilon$ but we shall for
simplicity drop temporarily this dependence in the notation. Now directly from
(\ref{ZetaExpansion}) we have:%
\[
\gamma_{n}=\left.  (-1)^{n}\frac{d^{n}}{ds^{n}}f(s)\right\vert _{s-1}%
\]
Then, after some elementary calculations, we get the main result of the
present paper:%
\begin{equation}
\fbox{$\gamma _{n}=\frac{n!}{\varepsilon ^{n}}\sum\limits_{k=n}^\infty
\frac{(-1)^{k}}{k!}\alpha _{k}S_{k}^{(n)}$}  \label{gamma gen}
\end{equation}
where $S_{k}^{(n)}$ are signed Stirling numbers of the first kind. Note that
in the literature there are different conventions concerning denotation and
indices of Stirling numbers which may be confusing. Here, following
\cite{Abramowitz}, we shall adopt the following convention involving the
Stirling numbers and the Pochhammer symbol:%
\[
\left(  x\right)  _{k}\equiv\frac{\Gamma(k+x)}{\Gamma(x)}=%
{\displaystyle\prod\limits_{i=0}^{k-1}}
(x+i)=(-1)^{k}%
{\displaystyle\sum\limits_{i=0}^{k}}
(-1)^{i}S_{k}^{(i)}x^{i}=%
{\displaystyle\sum\limits_{i=0}^{k}}
\left\vert S_{k}^{(i)}\right\vert x^{i}%
\]

Denoting%
\[
\beta_{nk}\equiv\frac{n!}{k!}\frac{S_{k}^{(n)}}{\varepsilon^{n}}%
\]
we can rewrite (\ref{gamma gen}) as formally an infinite matrix product%
\begin{equation}
\gamma_{n}=%
{\displaystyle\sum\limits_{k=n}^{\infty}}
\beta_{nk}\;\alpha_{k}\label{gamma gen 1}%
\end{equation}
The summation over $k$ starts from $n$ since $\beta_{nk}\equiv0$ for
$k<n$.\ Precision of $\alpha_{1}$ is equal to precision of precomputed values
of $f(s)$ given by (\ref{regularized}) in equidistant nodes. When $k$ grows
the precision of consecutive $\alpha_{k}$ almost linearly tends to zero. Thus
there always exists certain cut-off value of $k=k_{0}$. Therefore the
summation in (\ref{gamma gen 1}) should be performed to this value:%
\begin{equation}
\gamma_{n}=%
{\displaystyle\sum\limits_{k=n}^{k_{0}}}
\beta_{nk}\,\alpha_{k}\label{gamma gen 2}%
\end{equation}
(Adding more terms is inessential. In other words, one cannot compute
$\gamma_{n}$ for $n>k_{0}$. In order to perform this one should increase
precision of the precomputed values of $f(s)$ which in turn would
proportionally increase $k_{0}$, see Figure 9 for detailed description.)

As pointed earlier $\varepsilon$ need not to be small, however, choosing
smaller $\varepsilon$ greatly accelerates convergence of the series. Yet, it
turns out that smaller $\varepsilon$ implies smaller $k_{0}$. What is really
important: all significant digits of $\gamma_{n}$ obtained from the finite sum
(\ref{gamma gen 2}) are correct.%

\begin{figure}[ptb]%
\centering
\includegraphics[height=0.4\textheight]{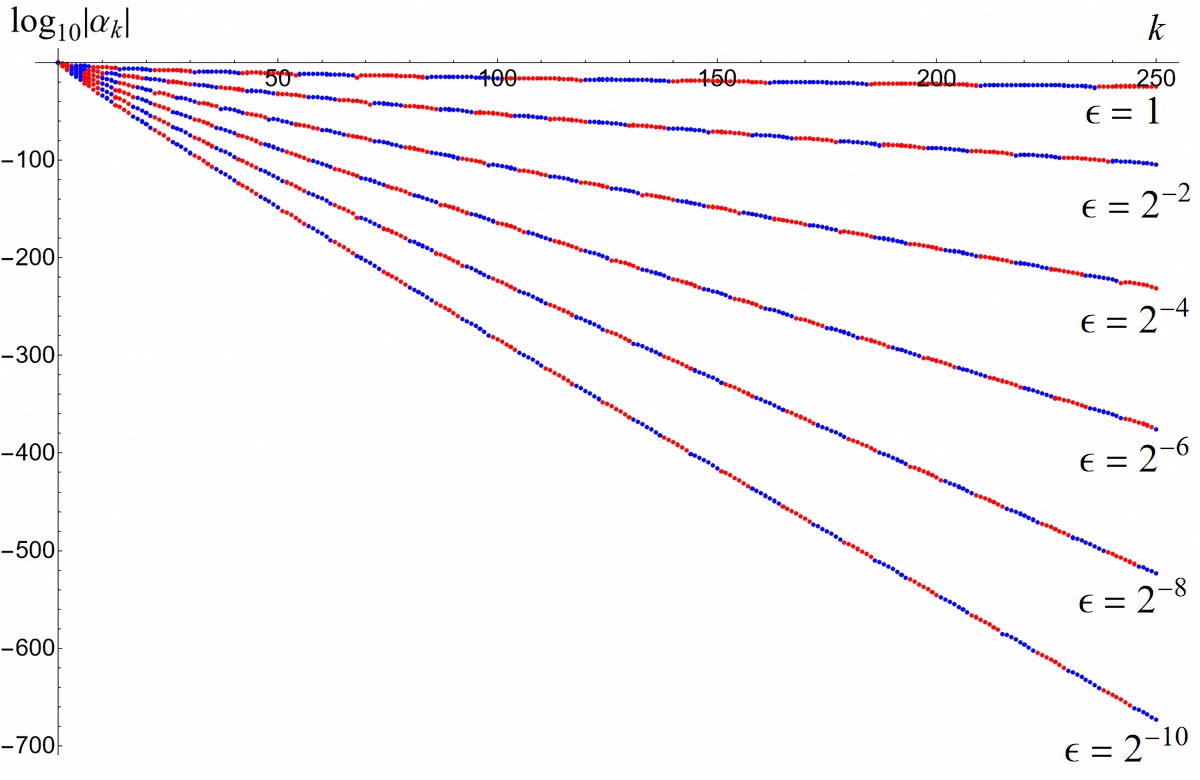}%
\caption{Behavior of coefficients $\alpha_{k}$ given by (\ref{Ak general}) for
different choices of the parameter $\varepsilon$.}%
\end{figure}
%

\begin{figure}[ptb]%
\centering
\includegraphics[height=0.4\textheight]{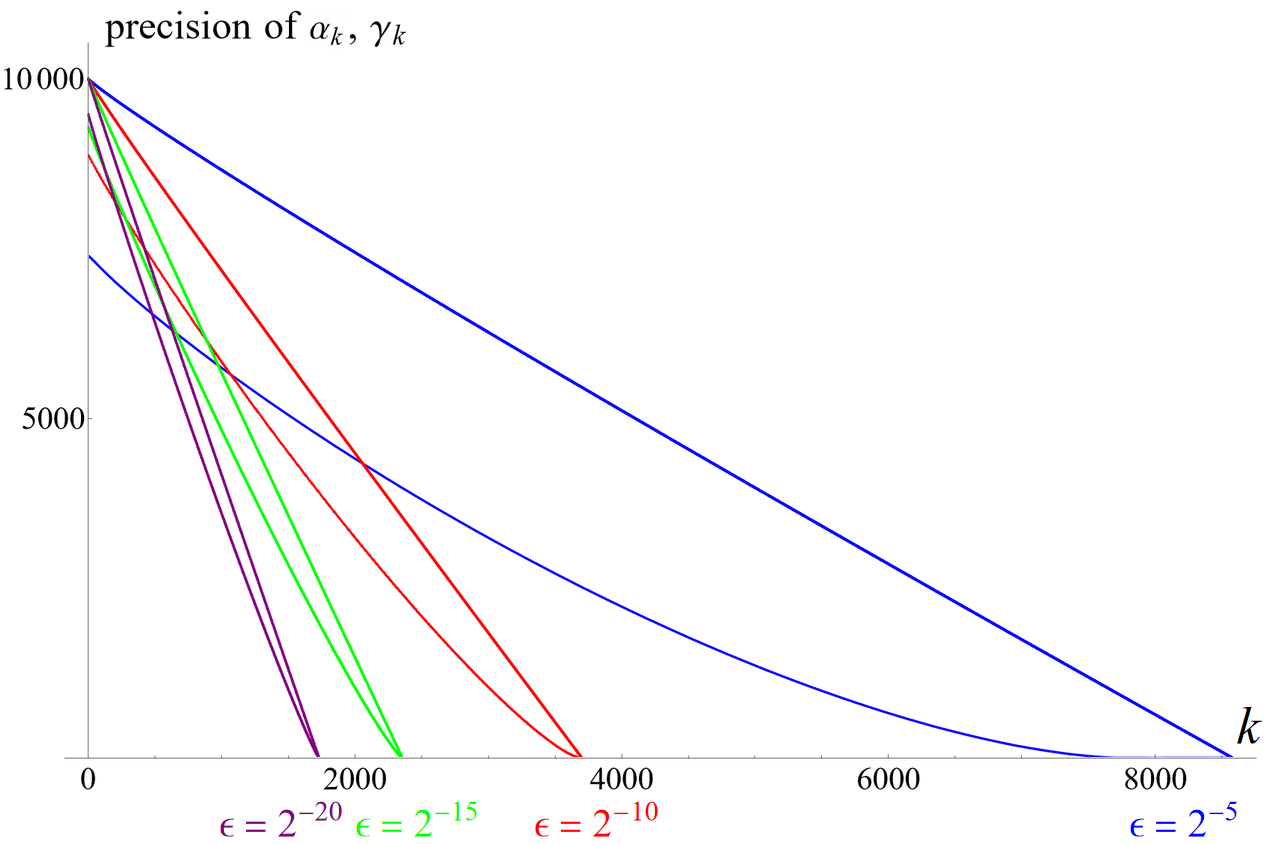}%
\caption{Plots of precision of $a_{k}$ versus $k$ for precision of the
precomputed equidistant zeta values equal to 10000 significant digits and for
four values of the auxiliary parameter $\varepsilon$: blue -- $\varepsilon
^{-5}$, red -- $\varepsilon^{-10}$, green -- $\varepsilon^{-15}$, purple --
$\varepsilon^{-20}$. For each color the upper, nearly straight line segment
corresponds to precision of $\alpha_{k}$, and piece of a curved line of the
same color corresponds to the precision of $\gamma_{k}$. Since precision of
$\alpha_{k}$ diminish with growing $k$ therefore for any given value of
$\varepsilon$ there always exist specific, unambiguous value of index $k_{0}$
such that for all $k>k_{0} $ precisions of all $\alpha_{k}$ are numerically
zero. Hence in the formula (\ref{gamma gen 2}) it is enough to sum only up to
this value.}%
\end{figure}

Of course, $\gamma_{n}$ eventually does not depend on particular choice of
$\varepsilon$, as expected, although $\alpha_{k}$ as well as the rate of
convergence of (\ref{gamma gen})--(\ref{gamma gen 2}) does. In fact series
(\ref{gamma gen}) converges for any value of $\varepsilon>0$\ but the rate of
convergence becomes terribly small for $\varepsilon\sim1$. On the other hand,
the smaller $\varepsilon$ the faster the rate of convergence. However, since
$\alpha_{k}$ also depends on $\varepsilon$, choosing smaller value for
$\varepsilon$ requires higher precision of precalculated values of $f(s)$
which in turn may be very time consuming. Hence, an appropriate compromise in
choosing $\varepsilon$ is needed.

Formula (\ref{gamma gen 2}) is particularly well-suited for numerical
calculations. Typically the algorithm has three simple steps:\bigskip

1. Tabulating function (\ref{regularized}) for equidistant arguments
$1+j\varepsilon$, i.e. $f(1+j\varepsilon),j=0,1,2,...$ This is most time
consuming and requires appropriate choice of parameter $\varepsilon$. (In our
case, we have chosen the value $\varepsilon=2^{-10}$.) The most convenient for
these calculations seems small but extremely efficient program PARI/GP which
has implemented particularly optimal zeta procedure. We used Cyfronet
Prometheus computer where calculating single value of $f(s)$ with 80000
significant digits requires about 10-15 minutes each. Since this procedure may
easily be parallelized therefore in order to compute more than 30000 values of
$f$ we started several dozen independent routines (each calculating a few
thousands values of $f$).

2. Calculating $\alpha_{k}$ using (\ref{Ak general}) and the precomputed
values of $f$.

3. Calculating Stieltjes constants using (\ref{gamma gen}).

(Contrary to the above steps 2. and 3., step 1 requires a powerful computer,
whereas steps 2 and 3 can be quickly performed on a typical PC.) It should be
emphasized that with the $\alpha_{k}$ coefficients properly calculated,
obtaining $\gamma_{n}$ requires only a dozen or so minutes on a very modest PC
machine. One property of the result (\ref{gamma gen}) should again be stressed
out: all digits of $\gamma_{n}$ obtained from (\ref{gamma gen 2}) are
significant and reliable.

Step 1 was achieved using the following PARI code:\bigskip%

$\backslash$%
g4%

$\backslash$%
p 80000;

default(parisizemax, 1000000000)

allocatemem(1000000000)

eps=2\symbol{94}-10;

f(s) = if ( s - 1, zeta(s) - 1/(s - 1) , Euler );

for( j = 0, 32000, write( "zeta.dat", "\{", 1+j*eps , "," , f(1+j*eps), "\},"));

\section{Appendix: Struggling with certain PARI bug}

Common experience shows that there are no computer programs, especially larger
ones, which -- in certain specific and usually unpredictable situations --
would not exhibit misbehavior. Computer program errors, according to the old
tradition called "bugs", are usually an integral part of each program. Of
course, program developers, or rather the large development teams that create
them, make every effort to ensure that their products are error-free. However,
it is virtually impossible to remove them completely. In addition,
professional computer programs are constantly developed and expanded,
sometimes over many years, and new functions are added in subsequent versions,
often at the explicit request of users. In this way, while previous bugs are
removed, new bugs are inevitably added, although, of course, this happens
unknowingly. The key role here is played by the fruitful cooperation of
program users with their developers: numerous users scattered all over the
world, solving their own specific problems, at the same time intensively test
programs and provide their developers with relevant information about
undesirable behavior of their products.

Certain sentence from a letter from PARI/GP user is particularly significant:
"I hope you and your collaborators will be able to eliminate the bugs [...] in
the forthcoming (final?) release of PARI" (April 1997). More than a quarter of
a century has passed since then. PARI is growing, has a faithful group of
users (mainly mathematicians dealing with the number theory), new functions
and procedures are added, but the list of "bugs" does not decrease at all. It
is instructive to look at the page%
\[
\text{https://pari.math.u-bordeaux.fr/cgi-bin/pkgreport.cgi?pkg=pari}%
\]
illustrating the intense and fruitful interaction of PARI users with its
creators. For someone unfamiliar with the essence of computer programs, the
sentence from the above-quoted letter sounds like the proverbial "wishful
thinking". It is naive and unrealistic, although it was sent in good faith.
The aforementioned "final version" of the program is an unattainable goal to
which one can, at best, "approach asymptotically". It is also worth adding
that this sentence was rightly placed on the PARI website with a meaningful
title: "Fun!".

When testing the algorithm described in this article, we came across a
surprising error in the numerical computation of the fundamental Riemann zeta
function, which is built into PARI. As mentioned earlier, the presented
algorithm requires "input" zeta values of great precision; in our case, we
chose $80000$ significant digits. It was a kind of compromise between
relatively high precision and reasonable computation time (several weeks on
many cores of the Prometheus supercomputer in Cyfronet in Krak\'{o}w).
Probably no one has methodically tested PARI for calculations of the Riemann
zeta with such great precision before.

The choice of PARI -- a small (in the command-line version only about 12 MB)
dedicated to calculations in number theory as a tool to obtain the value of
the zeta function -- resulted from the high speed of calculations: several
times greater than, for example, \textit{Mathematica} (size of installation
files over 4 GB). Unexpectedly tt turned out that the result file of the
necessary numerical values in the form of the array $\{1+j\varepsilon
,f(1+j\varepsilon)\}$ contained incorrect digits in the range of index $j$
from 11201 to 12401.

Of course, finding those digits that were wrong among more than 2.5 billion
digits was quite a challenge. It was a very tedious and frustrating job, like
looking for the proverbial needle in a haystack. But it was even more
challenging to figure out the very cause of this error. In the first case,
some properties of the $\alpha_{k}$ coefficients proved to be helpful. In the
second case, professional help of the PARI program developers turned out to be indispensable.

The first sign of the presence of these erroneous digits was that the
coefficients $\alpha_{k}$ calculated from these values, instead of rapidly
(exponentially) decrease to zero with the increase of the index $k$, changed
drastically its behavior, see Fig. 10.%
\begin{figure}[ptb]%
\centering
\includegraphics[height=0.35\textheight]{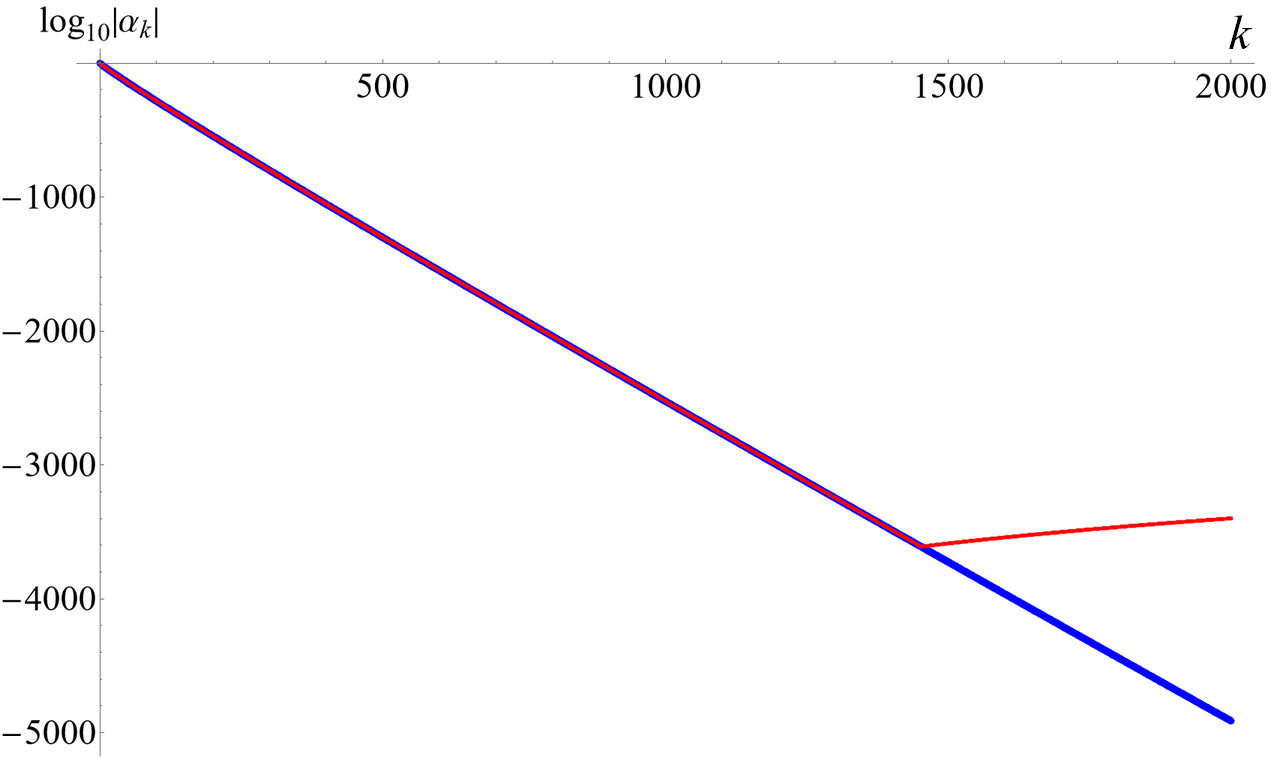}%
\caption{Coefficients $\alpha_{k}$ given by (\ref{Ak general}) are extremely
sensitive to even one wrong digit in the calculated value of the zeta
function, even at a very distant place of its decimal expansion. The figure
illustrates a sudden change in the behavior of the $\alpha_{k}$ coefficients
(red line) when in the correct value of the regularized zeta function
$f(1+\varepsilon j)$ for $j=1000$ only single digit is replaced with another
one that differs from the correct one just by $1$. The replaced digit might be
in a very remote significant place (in this case it was on position 4000 after
the decimal point), and yet $\alpha_{k}$ would "feel" and reveal that change
anyway.}%
\end{figure}

Since the (regularized) zeta function, however complicated and mysterious, is
a regular function, the successive finite differences of equidistant values of
this function from the above-mentioned table $f(1+j\varepsilon)$\ should lie
on a smooth curve. The tests performed with the use of the Mathematica
procedure \textbf{Differences} that calculates successive finite differences
revealed that for the above-mentioned values of index $j$ and with the order
of these differences about $400$, disturbing oscillations appeared instead of
a sequence of points lying along a smooth curve, see Fig. 11.%

\begin{figure}[ptb]%
\centering
\includegraphics[height=0.35\textheight]{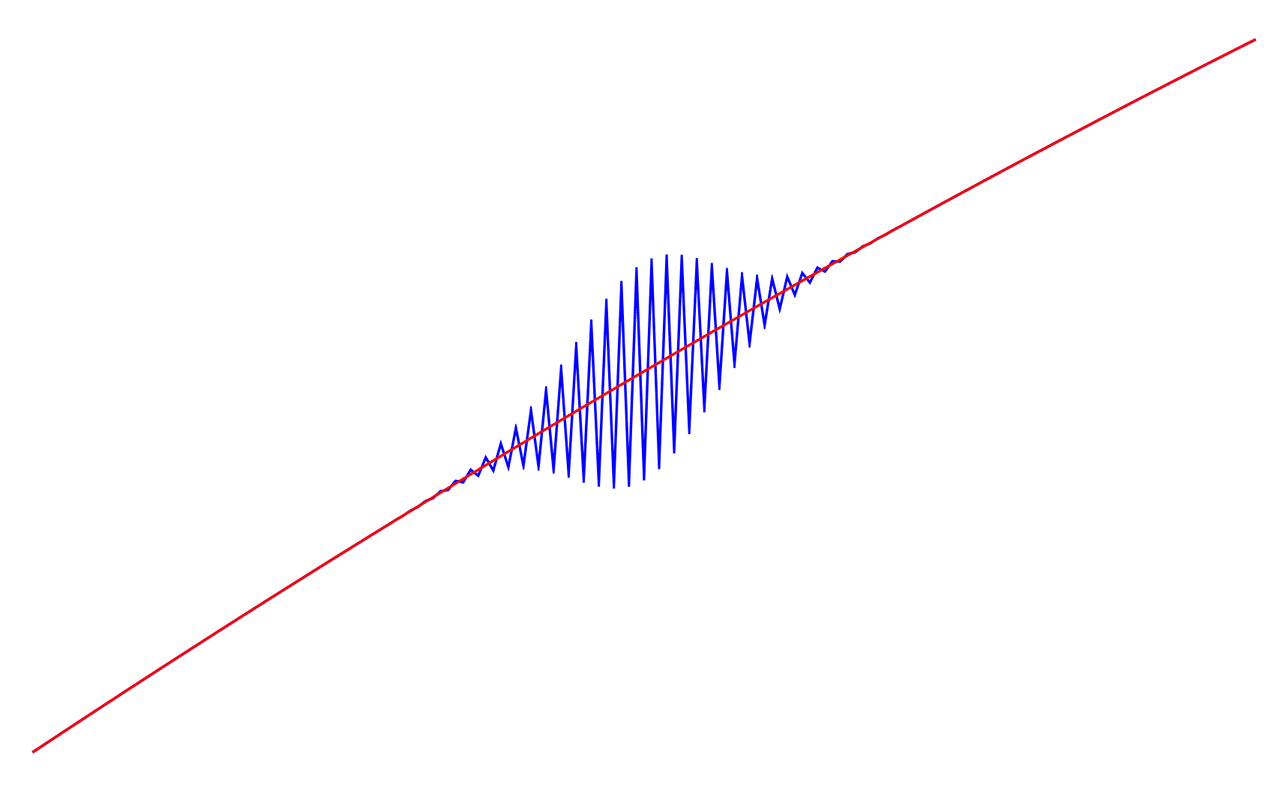}%
\caption{Wolfram Mathematica procedure \textbf{Differences[\textit{list}%
,\textit{n}], }which gives the $n^{\text{th}}$ differences of a given list,
when applied to the list of "contaminated" data of the function
$f(1+j\varepsilon)$ reveals that some digits are wrong. When the integer
parameter $n$ is sufficiently large then even a single wrong digit differing
from the proper one by unity produces oscillations instead of a smooth
distribution of points. In this case one had to use $n=400$ to reveal the
error. (In the above graph, both axes have been removed as they are irrelevant
to demonstrate the effect described.)}%
\end{figure}

Intensive and very tedious tests, requiring great patience, time and computer
resources, lasted for several weeks and were carried out with professional and
very kind cooperation of employees of the Cyfronet Computer Center in
Krak\'{o}w (administrators of the Prometheus supercomputer). In order to
eliminate the potential causes of generating wrong digits, we tested newer and
newer development versions of PARI released daily. We used two different
compilers (Intel icc 19.1.1.217 and GNU gcc version 4.8.5 20150623). We have
compiled PARI in serial and parallel version (threading engine: pthread, mpi,
single). Additionally, for the parallel version, we also ran single-core jobs
to rule out the PARI "parfor" command as a possible source of the problem. We
used different operating systems (Linux and Windows 10), different versions of
Linux cores (x86-64, x86-64 / GMP-6.2.1, x86-64 / GMP-6.0.0) and different
types of processors (Intel and AMD). We compared the obtained numbers with the
results obtained with Wolfram Mathematica on PCs with AMD and Intel processors
(these calculations took several times longer than with PARI). We additionally
performed a series of calculations for precision from 30000 to 90000 in 10000
steps and from 71000 to 89000 in 1000 steps. It turned out that the wrong
numbers appeared at 174 decimal places only for the precision of 74000 and 80000.

The results of these tests were successively (from August 2021 to December
2021) delivered to the authors of the PARI program, who made appropriate
corrections in the program code. (Incidentally, the first such correction did
not remove the error; it appeared again but in a different range of index $j$,
and even worse, i.e. for more significant digits of the Riemann zeta function...)

In the end, it turned out that the cause of this error was simply in the
PARI-implemented procedure for computing the value of the zeta function which
uses the classical Euler-Maclaurin algorithm. Specifically, the values of the
Bernoulli numbers required to compute the zeta function were rounded
unnecessarily. It was due to double roundings occurring when caching Bernoulli
numbers, because of too frequent precision reductions. This bug did not affect
the low precision computations, but was particularly bothersome with the
algorithm described here. More details can be found here:%
\[
\text{https://pari.math.u-bordeaux.fr/cgi-bin/bugreport.cgi?bug=2311}%
\]

It should be emphasized that when computing the zeta function, PARI first
computes and tabulates the appropriate Bernoulli numbers, according to the
Euler-Maclaurin formula. During this stage of the calculations, which --
depending on the precision set at the beginning -- sometimes takes several
hours, the results do not start to appear, and the program pretends that it
has "hung".

The revised version of the program finally appeared at the end of December
2021. From that moment, having the necessary and reliable numerical data, i.e.
the high precision (regularized) zeta values, we were able to return to purely
mathematical problems and continue the main project of calculating the
Stieltjes constants.

Finally, it should be emphasized once again that the main advantage of the
presented here algorithm for calculating important Stieltjes constants is its
mathematical simplicity and numerical efficiency. Moreover, it can be a
convenient starting point for deriving certain new asymptotic expansion for
these constants, both more accurate and simpler than several expansions known
in the literature. This will be the topic of another publication
\cite{Maslanka 3}.

\begin{acknowledgement}
This research was supported in part by PL-Grid Infrastructure and was
performed on the Prometheus multicore supercomputer, Cyfronet, Krak\'{o}w. The
kind and patient willingness to help and professional technical advice from
the members of the Cyfronet team (Maciej Czuchry and Andrzej Dorobisz from
High Performance Computing Software Department) when starting the project were
very helpful. Finally, the authors of the Pari/GP computer algebra system
(Bill Allombert and Karim Belabas) put a lot of effort into correcting the
persistent error in the calculation of the Riemann zeta function for the great
accuracies that were crucial to the presented algorithm.
\end{acknowledgement}

\end{document}